\documentclass{amsart}
\usepackage{amssymb}
\newtheorem{theorem}{Theorem}[section]
\newtheorem{conjecture}{Conjecture}
\newtheorem{corollary}[theorem]{Corollary}
\newtheorem{example}[theorem]{Example}
\theoremstyle{remark}
\newtheorem*{acknowledgements}{Acknowledgements}
\numberwithin{equation}{section}

\begin{document}
\title[Orthogonal harmonic analysis and scaling of fractal measures]{Orthogonal harmonic analysis and scaling of fractal measures\linebreak Analyse harmonique
orthogonale des mesures fractales avec structure d'\'echelle}
\author{\textsc{Palle E.T. Jorgensen} and \textsc{Steen Pedersen}}
\maketitle

\noindent\raisebox{150pt}[0pt][0pt]{Analyse math\'ematique/\it
Mathematical Analysis}

\section*{Abstract}

We show that certain iteration systems lead to fractal measures admitting
exact orthogonal harmonic analysis.

\section*{R\'{e}sum\'{e}}

On montre que certains syst\`{e}mes it\'{e}ratifs conduisent aux mesures
fractales qui admettent une analyse harmonique orthogonale exacte.

\section*{Version fran\c{c}aise abr\'{e}g\'{e}e}

Soit $\Omega$ un sous-ensemble mesurable pour la mesure de Lebesgue $m$ sur
l'espace Euclidien $\mathbb{R}^{d}$, $d\geq1$. Soit $L^{2}\left(  m_{\Omega
}\right)  $ l'espace de Hilbert des fonctions de carr\'{e} $m_{\Omega}%
$-int\'{e}grable par rapport au produit scalaire $\langle f\mid g\rangle
:=\int\overline{f(x)}\,g(x)\,dm_{\Omega}(x)$, o\`{u} $m_{\Omega}$ est la
restriction de la mesure de Lebesgue \`{a} $\Omega$. Le probl\`{e}me de savoir
pour quels sous-ensembles de mesure finie $\Omega$ il existe une base
orthogonale $\{e_{\lambda}(x):=\exp(i2\pi\lambda\cdot x):\lambda\in\Lambda\}$
dans $L^{2}\left(  m_{\Omega}\right)  $ a \'{e}t\'{e} soulev\'{e} par I.E.
Segal (1957) et, dans l'article \cite{Fugl74}, de B. Fuglede, et a
\'{e}t\'{e} \'{e}tudi\'{e} dans
\cite{Fugl74,Jorg82,JoPe92,LaWa97,Pede87,Pede96}. Il est bien connu (cf.
\cite{Fugl74,Pede87}) qu'un sous-ensemble ouvert, connexe, et de mesure finie
$\Omega$ de $\mathbb{R}^{d}$ admet une base orthogonale si et seulement si la
famille des d\'{e}riv\'{e}es partielles $-i\frac{\partial\,}{\partial x_{j}}$
qui op\`{e}rent sur $C_{c}^{\infty}(\Omega)$, l'espace des fonctions lisses et
\`{a} support compact dans $\Omega$, admet par extension une famille
d'op\'{e}rateurs hermitiens $H_{j}$, $1\leq j\leq d$, fortement deux-\`{a}%
-deux commutatifs (dans le sens qu'ils ont des r\'{e}solutions spectrales commutatives).

Quand $\Omega=\left[  0,1\right]  ^{d}$ est un cube dans $\mathbb{R}^{d}$, la
classe de tous ces op\'{e}rateurs a \'{e}t\'{e} d\'{e}couverte dans
\cite{JoPe97}, avec des r\'{e}sultats particuli\`{e}rement explicites pour
$d\leq3$. Des domaines $\Omega$ admettant de tels op\'{e}rateurs d'extension
mais qui satisfont seulement \`{a} une forme plus faible de commutativit\'{e}
ont \'{e}t\'{e} \'{e}tudi\'{e}s auparavant par J.~Friedrich \cite{Frie87} dans
le cas $d=2$.

On \'{e}tudie une mesure autosimilaire $\mu$ avec support contenu dans
l'intervalle $\left[  0,1\right]  $ et telle que le sous-espace vectoriel
engendr\'{e} par l'ensemble des fonctions analytiques $\left\{  e^{i2\pi
nx}:n=0,1,2,\dots\right\}  $ soit dense dans $L^{2}\left(  \mu\right)  $. On
identifie selon la dimension fractale de $\mu$ les sous-ensembles
$P\subset\mathbb{N}_{0}=\left\{  0,1,2,\dots\right\}  $ tels que les fonctions
$\left\{  e_{n}:=e^{i2\pi nx}:n\in P\right\}  $ constituent une base
orthogonale de $L^{2}\left(  \mu\right)  $. On donne aussi en dimension plus
grande une construction affine qui conduit aux mesures autosimilaires $\mu$
ayant leurs supports dans $\mathbb{R}^{d}$. Celle-ci est obtenue \`{a} partir
d'une matrice expansive d'ordre $d$ et d'un ensemble fini de vecteurs de
translation. En plus, pour que l'espace $L^{2}\left(  \mu\right)  $
correspondant ait une base orthogonale de fonctions exponentielles
$e^{i2\pi\lambda\cdot x}$ ayant des vecteurs $\lambda$ dans $\mathbb{R}^{d}$
comme ensemble d'indices, il faut que certaines conditions g\'{e}om\'{e}triques
(qui ont des rapports avec les op\'{e}rateurs de transfert de Ruelle) sur le
syst\`{e}me affine soient remplies.

On cite et discute ci-dessous quelques conjectures concernant les mesures qui
admettent une analyse harmonique orthogonale exacte.

Soit $\Omega$ un sous-ensemble mesurable de mesure finie pour la mesure de
Lebesgue sur $\mathbb{R}^{d}$. S'il existe un ensemble d'indices $\Lambda$ tel
que $\left\{  e_{\lambda}:\lambda\in\Lambda\right\}  $ est une base
orthogonale de $L^{2}\left(  \Omega\right)  $, alors on dit que $\Omega$ est un
\emph{ensemble spectral,} $\Lambda$ est le \emph{spectre,} et $\left(
\Omega,\Lambda\right)  $ est une \emph{paire spectrale.} S'il existe un
ensemble $T$ tel que, modulo des ensembles negligeables, la famille $\left\{
\Omega+t:t\in T\right\}  $ est une partition de $\mathbb{R}^{d}$, alors on dit
que $\Omega$ est un \emph{pav\'{e}} et que $T$ est un \emph{ensemble de
pavage. }On a alors la conjecture suivante de Fuglede.

\begin{conjecture}
[Conjecture de Fuglede \cite{Fugl74}]\label{Con1.1}Soit $\Omega$ un ensemble
de mesure finie et positive. Alors $\Omega$ est un ensemble spectral si et
seulement si $\Omega$ est un pav\'{e}.
\end{conjecture}

Cette conjecture reste ouverte dans les deux sens m\^{e}me pour le cas $d=1$.

Un ensemble $S$ qui est l'image d'un ensemble de la forme $\mathbb{Z}^{d}+A$
pour certain ensemble fini $A$ tel que $(A-A)\cap\mathbb{Z}^{d}=\varnothing$
est dit \emph{p\'{e}riodique}, ou un \emph{treillis avec une base.}

\begin{conjecture}
[Conjecture sur les ensembles spectraux p\'{e}riodiques \cite{Pede97}%
]\label{Con1.2}Soit $\Omega$ un ensemble de mesure finie et positive. Alors
$\Omega$ est un ensemble spectral admettant un spectre p\'{e}riodique si et
seulement si $\Omega$ est un pav\'{e} admettant un ensemble de pavage p\'{e}riodique.
\end{conjecture}

Il a \'{e}t\'{e} d\'{e}montr\'{e} (cf.\ \cite{Fugl74,Jorg82,Pede87}) que
$\Omega$ est un ensemble spectral ayant pour spectre $\mathbb{Z}^{d}$, si et
seulement si, $\Omega$ est un pav\'{e} avec l'ensemble de pavage
$\mathbb{Z}^{d}$. Les articles \cite{LaWa97} et \cite{Pede96} r\'{e}duisent la
conjecture sur les ensembles spectraux p\'{e}riodiques \`{a} certaines
questions concernant des sous-ensembles finis du treillis $\mathbb{Z}^{d}$.
Pour $d=1$ des progr\`{e}s ont \'{e}t\'{e} achev\'{e}s vers la r\'{e}solution
de ces probl\`{e}mes dans \cite{PeWa97}. Ces r\'{e}sultats supportent le point
de vue selon lequel certaines classes sp\'{e}cifiques d'ensembles spectraux
correspondent \`{a} certaines classes d'ensembles de pavage. D'ailleurs ceci
est confirm\'{e} par quelques r\'{e}sultats dans \cite{JoPe97}, o\`{u} nous
montrons que tout ensemble spectral p\'{e}riodique du cube $\Omega=[0,1]^{d}$
est aussi un ensemble de pavage pour le cube.

\section{Introduction}

\label{S:introduction}

We consider, in this note and in \cite{JoPe97} and \cite{JoPe97b},
generalized spectral transforms for a certain Fourier duality in
$\mathbb{R}^{d}$. Our results are motivated by considerations of the transform%
\[
\xi\longmapsto\int_{\Omega}e^{-i\xi\cdot x}f\left(  x\right)  \,dx,\quad\xi
\in\mathbb{R}^{d},\;f\in L^{2}\left(  \Omega\right)  ,
\]
for a given measurable subset $\Omega\subset\mathbb{R}^{d}$ of finite Lebesgue
measure. Instead we consider pairs of measures $\left(  \mu,\nu\right)  $ on
$\mathbb{R}^{d}$ such that the following generalized transform,%
\begin{equation}
F_{\mu}f\colon\lambda\longmapsto\int e^{-i2\pi\lambda\cdot x}f\left(
x\right)  \,d\mu\left(  x\right)  ,\label{eq1.1}%
\end{equation}
induces an isometric isomorphism of $L^{2}\left(  \mu\right)  $ onto
$L^{2}\left(  \nu\right)  $, specifically making precise the following
unitarity:%
\[
\int\left|  f\left(  x\right)  \right|  ^{2}\,d\mu\left(  x\right)
=\int\left|  \left(  F_{\mu}f\right)  \left(  \lambda\right)  \right|
^{2}\,d\nu\left(  \lambda\right)  .
\]
When applied to the case when $\mu$ is a measure of compact support with
\emph{fractal} Hausdorff dimension, we identify some candidates for pairs
$\left(  \mu,\nu\right)  $, in concrete examples, when duality does hold.

\section{Pairs of Measures}

\label{S:measures}

\subsection{New Pairs from Old Pairs\label{SS:NewOld}}

Let $\mu$ and $\nu$ be Borel measures on $\mathbb{R}^{d}$. We say that
$(\mu,\nu)$ is a \emph{spectral pair} if the map $F_{\mu}$ from (\ref{eq1.1})
above, defined for $f\in L^{1}\cap L^{2}(\mu)$, extends by continuity to an
isometric isomorphism mapping $L^{2}(\mu)$ onto $L^{2}(\nu)$. It was shown in
\cite{Pede87} that if $\mu$ is the restriction of Lebesgue measure to a
connected set of infinite measure, then the result on extensions of the
directional derivatives, described above, remains valid.

\subsection{Which Measures are Possible?\label{SS:Possible}}

It turns out that the class of measures that can be part of a spectral pair is
fairly limited: specifically, we have

\begin{theorem}
[Uncertainty Principle]\label{Thm2.2}Suppose $(\mu,\nu)$ is a spectral pair,
$f\in L^{2}(\mu)$, $f\neq0$, and $A,B\subset\mathbb{R}^{d}$. If $\Vert
f-\chi_{A}f\Vert_{\mu}\leq\varepsilon$ and $\Vert Ff-\chi_{B}Ff\Vert_{\mu}%
\leq\delta$, then $(1-\varepsilon-\delta)^{2}\leq\mu(A)\nu(B)$.
\end{theorem}

\begin{theorem}
[Local Translation Invariance]\label{Thm2.3}Suppose $(\mu,\nu)$ is a spectral
pair, and $t\in\mathbb{R}^{d}$. If $\mathcal{O}$ and $\mathcal{O}+t$ are
subsets of the support of $\mu$, then $\mu(\mathcal{O})=\mu(\mathcal{O}+t)$.
\end{theorem}

Our work on generalized spectral pairs is motivated by M.N. Kolountzakis and
J.C. Lagarias who in \cite{KoLa96} discuss related tilings of the real line
$\mathbb{R}^{1}$ by a function.

The following result establishes a direct connection to the spectral pairs
mentioned above.

\begin{theorem}
\label{Thm2.4}Suppose $(\mu,\nu)$ is a spectral pair. If $\mu(\mathbb{R}%
^{d})<\infty$, then $\nu$ must be a counting measure with uniformly discrete support.
\end{theorem}

\section{Fractal Measures}

\label{S:fractals}

\subsection{Dual Iteration Systems\label{SS:Dual}}

Consider a triplet $(R,B,L)$ such that $R$ is an expansive $d\times d$ matrix
with real entries, and $B$ and $L$ are subsets of $\mathbb{R}^{d}$ such that
$N:=\#B=\#L$,
\begin{align}
R^{n}b\cdot l &  \in\mathbb{Z},\text{ for any }n\in\mathbb{N},\;b\in B,\;l\in
L,\label{Compatibility}\\
H_{B,L} &  :=N^{-1/2}\left(  e^{i2\pi b\cdot l}\right)  _{b\in B,\,l\in
L}\text{is a unitary }N\times N\text{ matrix}.\label{Hadamard}%
\end{align}
We introduce two dynamical systems, $\sigma_{b}(x):=R^{-1}x+b$ and $\tau
_{l}(x):=R^{\ast}x+l$, and the corresponding ``attractors'', $X_{\sigma
}:=\left\{  \sum_{k=0}^{\infty}R^{-k}b_{k}:b_{k}\in B\right\}  $ and
\begin{equation}
\mathcal{L}=X_{\tau}:=\left\{  \sum_{k=0}^{n}R^{\ast k}l_{k}:n\in
\mathbb{N},\;l_{k}\in L\right\}  .\label{eq:Ldef}%
\end{equation}
The set $X_{\rho}$ is then the support of the unique probability measure which
solves the equation
\begin{equation}
\mu=N^{-1}\sum_{b\in B}\mu\circ\sigma_{b}^{-1}.\label{eq:mu}%
\end{equation}
We show that, under certain geometric assumptions, the exponentials
$E(\mathcal{L})=\{e_{\lambda}:\lambda\in\mathcal{L}\}$ form an orthogonal
basis for $L^{2}(\mu)$. It follows from the assumptions on $(R,B,L)$ that
$E(\mathcal{L})$ is orthogonal; so the question is whether or not these
exponentials span all of $L^{2}(\mu)$. If we set
\begin{equation}
\chi_{B}(t):=N^{-1}\sum_{b\in B}e_{b}(t),
\end{equation}
then the expansiveness property of $R$ and (\ref{eq:mu}) imply an explicit
product formula for the Fourier transform of $\mu$,
\begin{equation}
\widehat{\mu}(t):=\int\overline{e_{t}(x)}\,d\mu(x)=\prod_{k=0}^{\infty}%
\chi_{B}(R^{\ast\,-k}t),\label{E:InfiniteProd}%
\end{equation}
the convergence being uniform on bounded subsets of $\mathbb{R}^{d}$. We
introduce the function
\begin{equation}
Q(t):=\sum_{\lambda\in\mathcal{L}}\left|  \widehat{\mu}(t-\lambda)\right|
^{2},\quad t\in\mathbb{R}^{d},
\end{equation}
and the Ruelle operator $C$ given by
\begin{equation}
\left(  Cq\right)  (t):=\sum_{l\in L}\left|  \chi_{B}(t-l)\right|  ^{2}%
q(\rho_{l}(t)),\label{eq:Cdef}%
\end{equation}
where $\rho_{l}(x):=R^{\ast\,-1}(x-l)$. Both $Q$ and the constant function
$\mathbf{1}$ are eigenfunctions for the Ruelle operator $C$ with eigenvalue
$1$, and the issue becomes one of multiplicity. The attractor
\begin{equation}
X_{\rho}:=\left\{  \sum_{k=0}^{\infty}-R^{\ast\,-k}l_{k}:l_{k}\in L\right\}
,\label{rho-recur-dom}%
\end{equation}
corresponding to the system $\left\{  \rho_{l}\right\}  $, will also be used below.

\subsection{Orthogonal Bases\label{SS:Bases}}

Let $H_{2}(\mathcal{L})$ denote the subspace of $L^{2}(\mu)$ spanned by the
orthonormal set $\{e_{\lambda}:\lambda\in\mathcal{L}\}$. Any $e_{t}$,
$t\in\mathbb{C}^{d}$ is in $L^{2}(\mu)$ so $H_{2}(\mathcal{L})$ is a subspace
of $L^{2}(\mu)$. We will show that $H_{2}(\mathcal{L})=L^{2}(\mu)$ for
specific systems $(R,B,L)$ satisfying (\ref{Compatibility})--(\ref{Hadamard}).

Let $Y$ denote the convex hull of the attractor $X_{\rho}$ given by
(\ref{rho-recur-dom}), and let $\left\|  q\right\|  _{\infty}:=\sup_{y\in
Y}\left|  q(y)\right|  $. We then introduce the following Lipschitz norm:%
\begin{equation}
\left\|  q\right\|  _{Y,\infty}:=\left\|  \left|  \nabla q\right|
_{2}\right\|  _{\infty}.\label{eq:Y-sup-norm}%
\end{equation}
We show that, if the operator norm of $C$, acting on a suitable set of smooth
functions, is less than one, then $\mu$ has the basis property.

\begin{theorem}
\label{T:Cbound}Let $(R,B,L)$ be a system in $\mathbb{R}^{d}$ satisfying
\textup{(\ref{Compatibility})--(\ref{Hadamard})}, $0\in L$. Let $C$ be the
operator given by \textup{(\ref{eq:Cdef})}, let $Y$ denote the convex hull of
the attractor $X_{\rho}$ given by \textup{(\ref{rho-recur-dom})}, and let
$\left\|  q\right\|  _{Y,\infty}$ be given by \textup{(\ref{eq:Y-sup-norm})}.
Supposing that $L$ spans $\mathbb{R}^{d}$, if there exists $\gamma<1$ such
that $\left\|  Cq\right\|  _{Y,\infty}\leq\gamma\left\|  q\right\|
_{Y,\infty}$ for all $q$ in a set of $C^{1}$-functions containing
$\mathbf{1}-Q$, then $H_{2}(\mathcal{L})=L^{2}(\mu)$.
\end{theorem}

The following result on Lipschitz estimates allows us to compute an explicit
and numerical operator norm bound $\gamma$ for $C$ in terms of the given data
$(R,B,L)$:

\begin{theorem}
\label{T:estimate}Let $(R,B,L)$ be a system in $\mathbb{R}^{d}$ satisfying
\textup{(\ref{Compatibility})--(\ref{Hadamard})}, $0\in L$. Let $C$ be the
operator given by \textup{(\ref{eq:Cdef})}, and let $Y$ denote the convex hull
of the attractor $X_{\rho}$ given by \textup{(\ref{rho-recur-dom})}. Let
$\left\|  q\right\|  _{Y,\infty}$ be given by \textup{(\ref{eq:Y-sup-norm})},
and
\[
\beta:=2\pi\operatorname{diam}(B)\vphantom{\max_{b,b^{\prime}\in B}}%
\smash{\max_{\substack{b,b^{\prime}\in B\\l\in L}}}\left\|  \sin
(2\pi(b-b^{\prime})(\,\cdot\,-l))\right\|  _{\infty}.
\]
Then we have the estimate
\[
\left\|  Cq\right\|  _{Y,\infty}\leq\left[  \left(  N-1\right)  ^{2}%
N^{-1}\beta\left\|  R^{-1}\right\|  _{\operatorname{op}}\max_{l\in L}\left|
l\right|  _{2}+\left\|  R^{-1}\right\|  _{\operatorname{hs}}\right]
\,\left\|  q\right\|  _{Y,\infty},
\]
for all $C^{1}$-functions $q$ such that $q(0)=0$. Here $\left\|  T\right\|
_{\operatorname{op}}$ is the operator norm, and $\left\|  T\right\|
_{\operatorname{hs}}$ is the Hilbert--Schmidt norm for a $d\times d$ matrix
$T$.
\end{theorem}

\begin{corollary}
\label{Cor3.3}Let $(R,B,L)$ satisfy \textup{(\ref{Compatibility}%
)--(\ref{Hadamard})}, and for $r\in\mathbb{N}$ let
\[
\mathcal{L}_{r}:=\left\{  \sum_{k=0}^{n}\left(  rR^{\ast}\right)  ^{k}%
l_{k}:n\in\mathbb{N},l_{k}\in L\right\}
\]
and let $\mu_{r}$ be the probability measure solving $\mu_{r}=N^{-1}\sum_{b\in
B}\mu_{r}\circ\sigma_{r,b}^{-1}$, where $\sigma_{r,b}(x):=\left(  rR\right)
^{-1}x+b$ \textup{(}i.e., a scaled version of \textup{(\ref{eq:mu}))}. If $L$
spans $\mathbb{R}^{d}$ and $0\in L$, then, provided $r$ is sufficiently large,
it follows that $\{e_{\lambda}:\lambda\in\mathcal{L}_{r}\}$ is an orthonormal
basis for $L^{2}(\mu_{r})$.
\end{corollary}

R.S. Strichartz obtained an asymptotic and quite different harmonic analysis
for the class of measures considered in this paper: it was based instead on a
continuous transform (see \cite{Stri94} for a survey of Strichartz's work on
self-similarity in harmonic analysis).

\section{Applications}

\label{S:applications}

\subsection{Fractal Hardy Spaces\label{SS:Hardy}}

One way to construct systems $(R,B,L)$ satisfying (\ref{Compatibility}%
)--(\ref{Hadamard}) is to pick $R$, $B$ and $L$ so that
\begin{equation}
R\in M_{d}(\mathbb{Z}),\quad RB\subset\mathbb{Z}^{d},\quad L\subset
\mathbb{Z}^{d}. \label{Compatibility2}%
\end{equation}
In fact (\ref{Compatibility2}) implies (\ref{Compatibility}) since
$R^{n}b\cdot l=Rb\cdot R^{\ast(n-1)}l$ for $n=1,2,3,\ldots$. The only
condition that is hard to satisfy is (\ref{Hadamard}). This last condition is
notoriously difficult: for example, it is not known which matrices with
entries in the unit circle satisfy (\ref{Hadamard}) for $N=7$; see, e.g.,
\cite{BjSa95,haag1} for some progress in the study of (\ref{Hadamard}). The
condition (\ref{Compatibility2}) is closely related to a condition use in the
study of certain multi-dimensional wavelets. Some results for systems
$(R,B,L)$ satisfying (\ref{Compatibility2}) and (\ref{Hadamard}) were
established in \cite{JoPe96a}.

If $R$ has non-negative integer entries, we will often end up with
$\{e_{\lambda}:\lambda\in\mathcal{L}\}$ being an orthonormal basis for
$L^{2}(\mu)$, and each element in $\mathcal{L}$ only having non-negative
coordinates. This is an interesting situation because the basis property leads
to the expansion $f=\sum_{\lambda\in\mathcal{L}}\left\langle e_{\lambda}\mid
f\right\rangle _{\mu}e_{\lambda}$, so setting $z_{j}:=e^{i2\pi x_{j}}$ we see
that $f(x)=\sum_{\lambda\in\mathcal{L}}\left\langle e_{\lambda}\mid
f\right\rangle _{\mu}z^{\lambda}$ for $f\in L^{2}(\mu)$, where $z^{\lambda
}:=\prod_{k=1}^{d}z_{k}^{\lambda_{k}}$. It follows that $f(x)$, $x\in
X_{\sigma}$, is the a.e.\ boundary value of a function analytic in the
polydisc $\{z\in\mathbb{C}^{d}:\left|  z_{j}\right|  <1\}$. Hence our
construction shows that many fractal $L^{2}$-spaces are Hardy spaces. This is
in sharp contrast to the Lebesgue spaces, for example, if $m_{[0,1]}$ is
Lebesgue measure restricted to the unit interval $[0,1]$.

\subsection{Examples\label{SS:Examples}}

Using Theorem \ref{T:estimate}, Theorem \ref{T:Cbound}, and equation
(\ref{E:InfiniteProd}), one can prove the following result.

\begin{theorem}
\label{T:DimOne} Suppose $d=1$, $N=2$, $B=\{0,a\}$, with $a\in\mathbb{R}%
\setminus\{0\}$, $R$ is an integer with $\left|  R\right|  \geq2$, and $\mu$
is given by \textup{(\ref{eq:mu})}. If $R$ is odd, then $L^{2}(\mu)$ does not
have a basis of exponentials for any $a\in\mathbb{R}\setminus\{0\}$. If $R$ is
even and $\left|  R\right|  \geq4$, then $L^{2}(\mu)$ has a basis of
exponentials for all $a\in\mathbb{R}\setminus\{0\}$.
\end{theorem}

Using convolution and Theorem \ref{T:DimOne} one can verify the following example.

\begin{example}
Let $\mu_{0}$ be the probability measure solving \textup{(\ref{eq:mu})} when
$R=4$ and $B=\{0,1/2\}$. Let $L=\{0,1\}$ and let $\mathcal{L}$ be given by
\textup{(\ref{eq:Ldef})}. Set $\Omega:=[0,1]+\mathcal{L}$, and define two
measures $\mu$ and $\nu$: $\mu(\Delta):=m(\Delta\cap\Omega)$, $\nu
(\Delta):=\sum_{k=0}^{\infty}\mu_{0}(\Delta+k)$. Then $(\mu,\nu)$ is a
spectral pair, and $\Omega$ is a tile with tiling set $-2\mathcal{L}$.
\end{example}

This is an example of a spectral set of infinite measure whose spectrum is not
periodic, and it takes us full circle, connecting back to the two extension
problems discussed in Conjectures \ref{Con1.1} and \ref{Con1.2} above.

\begin{acknowledgements}
We wish to thank Brian Treadway for the beautiful typesetting, and Tuong Ton
That for help with the French part. Partial support from the National Science
Foundation is gratefully acknowledged.
\end{acknowledgements}

\bigskip

\noindent \begin{minipage}[t]{0.5\textwidth}
\noindent Palle E.T. Jorgensen\newline Department of Mathematics\newline The
University of Iowa\newline Iowa City, IA
52242-1419\newline U.S.A.\newline
T\'el.\ :\ (1) 319.335.0782\newline
Fax\ :\ (1) 319.335.0627\newline
E.mail\ :\ \texttt{jorgen@math.uiowa.edu}\end{minipage}%
\begin{minipage}[t]{0.5\textwidth}
\noindent Steen Pedersen\newline Department of Mathematics\newline
Wright State University\newline
Dayton, OH 45435\newline
U.S.A.\newline
T\'el.\ :\ (1) 937.775.2432\newline
Fax\ :\ \newline
E.mail\ :\ \texttt{steen@math.wright.edu}\end{minipage}\bigskip

\noindent Please send proofs to/Veuillez envoyer les \'epreuves \`a\ :\ Prof. Jorgensen \textit{supra}.
\end{document}